\documentclass{siamart251104}
\usepackage{amsfonts}
\usepackage{hyperref}
\usepackage{cleveref}
\usepackage{stmaryrd}
\usepackage{graphicx}
\usepackage[svgnames]{xcolor}

\numberwithin{equation}{section}
\newcommand{\eps}{\varepsilon}
\newcommand{\jump}[1]{\llbracket #1 \rrbracket}
\newcommand{\average}[1]{\langle #1 \rangle}
\newcommand{\Real}{\mathrm{Re}}
\newcommand{\Image}{\mathrm{Im}}

\newcommand{\Add}[1]{#1}
\newcommand{\Erase}[1]{}

\newsiamthm{result}{Main Result}

\title{Dynamics of interfaces in the two-dimensional wave-pinning model}
\author{
    Shunsuke Kobayashi\thanks{Faculty of Engineering, University of Miyazaki, 1‑1 Gakuen Kibanadainishi, Miyazaki 889-2192, Japan}
    \and Koya Sakakibara\thanks{Faculty of Mathematics and Physics, Institute of Science and Engineering, Kanazawa University, Kakuma-machi, Kanazawa-shi, Ishikawa 920-1192, Japan}
    \and Taikei Uechi\thanks{Division of Mathematical and Physical Sciences, Graduate School of Natural Science and Technology, Kanazawa University, Kakuma-machi, Kanazawa-shi, Ishikawa 920-1192, Japan}
}

\headers{Dynamics of interfaces in the 2D WP model}{Shunsuke Kobayashi, Koya Sakakibara, and Taikei Uechi}

\begin{document}

\maketitle

\begin{abstract}
We study the mass-conserved reaction-diffusion system known as the wave-pinning model, which serves as a minimal framework for describing cell polarity.
In this model, the interplay between reaction kinetics and slow diffusion forms a sharp interface that partitions the domain into high- and low-concentration regions.
We perform a detailed asymptotic analysis and derive higher-order approximation equations governing the motion of this interface.
Our results show that on a fast timescale, the interface evolves via propagating front dynamics, whereas on a slow timescale, it evolves as an area-preserving mean curvature flow.
Furthermore, using the derived free boundary problem, we demonstrate that on a significantly slower timescale, an interface whose endpoints lie on the domain boundary drifts along the boundary toward regions of higher curvature.
In summary, our analysis reveals that the interface dynamics in the wave-pinning model exhibit a hierarchy of three distinct timescales: wave propagation on a fast timescale, curvature-driven area-preserving evolution on a slow timescale, and motion along the boundary on a significantly slower timescale.
\end{abstract}



\section{Introduction}
Cell polarity, characterized by the asymmetric distribution of cellular components, plays a fundamental role in a wide range of biological processes including cell migration, division, and morphogenesis \cite{bryant2008cells,butler2017planar,campanale2017development,wodarz2007cell}.
At the same time, its disruption is a hallmark of many pathological conditions such as cancer \cite{butler2017planar,lee2008cell,wodarz2007cell}.
Consequently, elucidating the mechanisms by which cells break symmetry remains a central challenge in modern biology.
In parallel with this biological inquiry, symmetry breaking has long been a subject of intensive mathematical research.
In this context, cell polarity has served as an intriguing case study, leading to the proposal of numerous mathematical models \cite{altschuler2008spontaneous,goehring2011polarization,goryachev2008dynamics,mori2008wave,otsuji2007mass}.
These studies have provided theoretical insights into how cells establish polarity and continue to drive research in mathematical biology (see, e.g., the reviews \cite{edelstein2013simple,goryachev2017many,jilkine2011comparison}).

In this paper, we study the wave-pinning model, which has emerged as a minimal and analytically tractable framework for explaining how cells achieve polarity \cite{mori2008wave,Mori2011}.
Biologically, this model captures the dynamics of Rho GTPases, which cycle between active membrane-bound and inactive cytosolic forms.
These dynamics are distilled into a mass-conserved reaction-diffusion system on a bounded domain $\Omega \subset \mathbb{R}^2$:
\begin{subequations}
    \label{wp}
    \begin{align}
        &\varepsilon\partial_t u = \varepsilon^2\Delta u + f(u,v) && \quad \text{in } \Omega\times (0,\infty), \label{wp-u}\\
        &\varepsilon\partial_t v = D\Delta v - f(u,v) && \quad \text{in } \Omega\times (0,\infty), \label{wp-v}\\
        &\partial_{\nu}u = \partial_{\nu}v = 0 && \quad \text{on } \partial\Omega\times (0,\infty),\label{wp-bc}
    \end{align}
\end{subequations}
where $\eps > 0$ is a sufficiently small parameter, $D > 0$ is a constant, and $\nu$ is the outward unit normal vector to the boundary $\partial\Omega$.
The functions $u$ and $v$ represent the concentrations of the active and inactive forms, respectively.
Detailed assumptions for this model are given in Section 2.
The generation of polarity in this framework hinges on three essential properties: (i) a large disparity in diffusion coefficients ($\eps^2 \ll D$), reflecting the rapid diffusion of cytosolic proteins; (ii) reaction kinetics $f(u,v)$ that are bistable with respect to $u$, allowing the coexistence of high and low states of $u$; and (iii) the conservation of total mass.
To understand the mechanism intuitively, imagine a local increase in the active form $u$. This initiates a wave of activation that propagates across the domain. However, since the total mass is conserved, the expansion of the high-$u$ region rapidly depletes the available inactive form $v$. As the level of $v$ drops below a critical level, the wave speed stalls, ``pinning'' the interface in place. This self-limiting growth is the essence of the wave-pinning mechanism.
To elucidate the detailed dynamics of this wave-pinning process, formal asymptotic analysis has been employed to characterize the front behavior as $\eps$ tends to $0$ \cite{Mori2011}.
Complementing these formal results, Gomez, Lam, and Mori rigorously demonstrated that the front solution to the shadow wave-pinning system converges to this sharp interface limit \cite{gomez2023front}.
Beyond these fundamental analyses, the model has been extended to incorporate various biological factors, such as source and sink terms \cite{verschueren2017model}, feedback from F-actin \cite{holmes2012regimes}, mechanical tension \cite{zmurchok2018coupling}, and bulk-surface coupling \cite{diegmiller2018spherical,giese2015influence,ratz2012turing}.
In addition, the model has been applied to studies of cell motility \cite{camley2013periodic,camley2017crawling}.

\begin{figure}[t]
    \centering
    \includegraphics[scale=0.37]{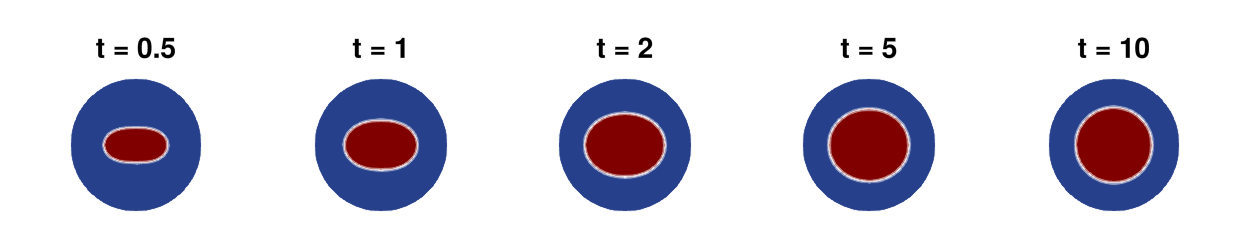}\label{fig:curvature-flow}
    \caption{Numerical simulation of \cref{wp} demonstrating the area-preserving, curvature-driven evolution on a slow timescale.}
\end{figure}

Here, we focus on the detailed dynamics of the interface after the initial front formation.
The interface equations derived in previous studies, such as \cite{Mori2011} and \cite{gomez2023front}, correspond to leading-order approximations that are independent of the interfacial mean curvature.
As noted by \cite{gomez2023front}, these approximations can exhibit finite-time singularities, underscoring the necessity for deriving more accurate, higher-order interface equations.
Incorporating such higher-order corrections is expected to regularize the dynamics and reveal that the interface evolves on a slow timescale analogously to volume-preserving mean curvature flow, as suggested by the literature \cite{gomez2023front,jilkine2009wave,mori2008wave}; indeed, our numerical simulations corroborate this behavior (Fig.~\ref{fig:curvature-flow}).
Similarly, for the bulk-surface coupled wave-pinning model, Miller et al. employed timescale separation techniques to demonstrate that the interface evolves as an area-preserving geodesic curvature flow \cite{miller2023generation}.
However, for the original wave-pinning model, a rigorous analysis of the slow dynamics is lacking.
Consequently, the derivation of a unified interface equation that remains valid across multiple timescales is an open problem.

\begin{figure}
    \centering
    \includegraphics[scale=0.45]{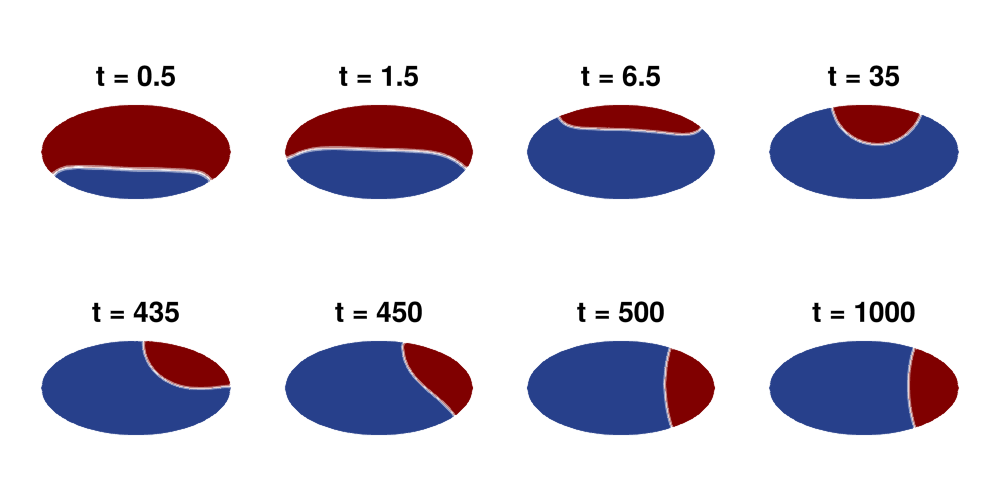}\label{fig:long-time-behavior}
    \caption{Numerical simulation of \cref{wp} on an elliptic domain. The dynamics exhibit a hierarchy of timescales: fast wave propagation, slow area-preserving curvature flow, and very slow drift along the boundary.}
\end{figure}

Through numerical simulations, another intriguing phenomenon regarding the interface motion of the wave-pinning model has been identified: the ``shape-sensing phenomenon'' \cite{maree2012cells,singh2022sensing,vanderlei2011computational} (Fig.~\ref{fig:long-time-behavior}).
In this process, high-activity domains drift very slowly along the boundary of the domain $\Omega$ toward regions of higher curvature.
Notably, the timescale of this phenomenon is distinct from, and typically slower than, the curvature flow regime described above.
Such geometry-dependent slow motion is not unique to the wave-pinning model; it has been observed in various reaction-diffusion systems, including the Allen-Cahn equation, suggesting it is a generic behavior driven by the interaction between the interface and the boundary.
This type of geometry-dependent slow evolution is classically referred to as ``slow dynamics'' or ``metastability'', which has been extensively studied in the mathematical analysis of reaction-diffusion systems \cite{ei1998slow,fusco1989slow}, and of related mass-conserving phase-separation models such as the Cahn-Hilliard equation and the Mullins-Sekerka flow \cite{Alikakos2000,alikakos2000motion,alikakos1998slow}.
Biologically, shape-sensing is a well-observed phenomenon; this is summarized in \cite{singh2022sensing}. In parallel, it is suggested to play a pivotal role in the phase-field model of cell motility proposed by \cite{camley2017crawling}.
In their model, the migration of the high-activity domain along the cell boundary is essential for driving rotational cell motion.
Therefore, providing a theoretical understanding of the internal wave-pinning dynamics---specifically, how the interface interacts with the domain geometry---would offer valuable insights.

Motivated by these considerations, we perform a derivation of more detailed interface equations in $\mathbb{R}^2$ using the method of matched asymptotic expansions.
While Jilkine pioneered the derivation of interface equations incorporating slow timescales \cite{jilkine2009wave}, we further develop this analysis to explicitly demonstrate the emergence of area-preserving mean curvature flow.
More precisely, we derive the following free boundary problem for the interface $\Gamma(t),\,t > 0$:
\begin{subequations}
    \label{fbp}
    \begin{align}
        &V_{n} = c(v_{0}) - \varepsilon[\kappa - \average{\kappa} + \beta(v_{0})(v_{1} - \average{v_1})] & &\text{on }\Gamma(t), &t \in (0,\infty), \label{fbp-ie}\\
        &u_{+}(v_{0})|\Omega_{+}(t)| + u_{-}(v_{0})|\Omega_{-}(t)| + v_{0}|\Omega| = M, & & &t \in (0,\infty), \label{fbp-mc0} \\
        &-\Delta v_{1} = \frac{(1 + u_{\pm}'(v_{0}))\jump{u_0}}{D \Lambda(v_{0})}c(v_0)L\Add{(t)} & &\text{in }\Omega_\pm(t), &t \in (0,\infty), \label{fbp-vp}\\
        &D\jump{\partial_{n} v_{1}} + c(v_{0})\jump{u_{0}} = 0 & &\text{on }\Gamma(t), &t \in (0,\infty), \label{fbp-jump}\\
        &\partial_{\nu}v_{1} = 0 & &\text{on }\partial\Omega, &t \in (0,\infty), \label{fbp-nf}\\
        &\int_{\Omega}(u_1 + v_1)dx = 0, \ u_1 = \frac{\partial_t u_0 - f_v(u_0, v_0)v_1}{f_u(u_0,v_0)}, & & &t \in (0,\infty). \label{fbp-mc1}
    \end{align}
\end{subequations}
Here, $\Omega_\pm(t)$ denote the subdomains where $u$ is close to the stable states $u_\pm(v_0)$, and we define $u_0 = u_\pm(v_0)$ in $\Omega_\pm(t)$.
$v_0$ and $v_1$ are the leading- and first-order terms of $v$.
$|\Omega_\pm(t)|$ denotes the area of these subdomains $\Omega_\pm(t)$.
The vector $n$ is the unit normal to the interface pointing from $\Omega_+(t)$ to $\Omega_-(t)$.
Accordingly, $V_n$ and $\kappa$ represent the normal velocity and the mean curvature of the interface, while $\average{\cdot}$ and $\jump{\cdot}$ denote the average and the jump (defined as the value in $\Omega_+$ minus that in $\Omega_-$) across the interface.
$c(v_0)$ is the propagating wave speed, $M$ the total mass, $D$ the diffusivity, and $L$ the interface length.
Finally, the coefficients $\beta(v_0)$ and $\Lambda(v_0)$ are defined in \cref{B(v_0)} and \cref{Lambda}, respectively.
$f_u$ and $f_v$ denote the partial derivatives of $f$ with respect to $u$ and $v$, respectively.
While we focus on the two-dimensional case, the same argument holds in general $N$ dimensions, with $\kappa$ replaced by $(N-1)\kappa$.
To the best of our knowledge, the derivation of the interface equations that approximate dynamics occurring across multiple timescales has not been established in the literature, with the notable exception of the wave-pinning model. We highlight this as a key mathematical novelty of the present work.

Furthermore, based on the derived free boundary problem, we investigate the behavior of an interface whose endpoints lie on the domain boundary.
We demonstrate that such an interface drifts along the boundary toward regions of higher curvature with a speed of $O(\varepsilon^2)$.
Specifically, assuming the interface can be approximated by a sufficiently small circular arc, we represent its center using the arc-length parametrization $X^0(s)$ of the boundary $\partial\Omega$.
We then show that the evolution of the position $s = s(t)$ is governed by the following approximation:
\begin{equation}
    \frac{d}{dt}s(t) = \frac{4\eps^2}{3\pi} K_{s}(s(t)),
\end{equation}
where $K_s$ denotes the derivative of the curvature $K$ of the boundary $\partial\Omega$ with respect to arc length.
In summary, the interface evolves as a propagating front on the $O(1)$ timescale, undergoes area-preserving mean curvature flow on the $O(\eps^{-1})$ timescale, and drifts along the boundary on the $O(\eps^{-2})$ timescale.
This hierarchy corresponds to a biological narrative of polarity establishment: First, the cell rapidly establishes a polarized front (fast timescale). Second, the polarized domain relaxes its shape to minimize interfacial tension while preserving its size (slow timescale). Finally, on a much longer timescale, the cell ``senses'' its geometric confinement, drifting toward regions of high curvature to stabilize its orientation (very slow timescale).

The organization of this paper is as follows.
In Section 2 we outline the model assumptions and derive the interface equations \cref{fbp}.
In Section 3 we investigate the slow motion of the interface along the domain boundary.
We conclude in Section 4 with some final remarks.

\section{Formal derivation of the interface equation}
In this section we present a formal derivation of the free boundary problem \cref{fbp} for the wave-pinning model \cref{wp}.
We perform the subsequent analysis under the following three assumptions for the reaction term $f$, adopting the formulation by \cite{Mori2011}:
\begin{description}
    \item[Bistability condition.] For any given $v$ in a certain range  $[v_\mathrm{min}, v_\mathrm{max}]$, the function $f(\cdot,v)$ is bistable.
This means that the equation $f(u,v) = 0$ for $u$ has three roots $u_-(v) < u_m(v) < u_+(v)$ for each $v \in [v_\mathrm{min}, v_\mathrm{max}]$, satisfying
\begin{equation}
    f_u(u_\pm(v),v) < 0, \quad f_u(u_m(v),v) > 0. \label{bistability-condition}
\end{equation}

    \item[Homogeneous stability condition.] The homogeneous states $(u,v) \equiv (u_\pm(v),v)$ for $v \in [v_\mathrm{min}, v_\mathrm{max}]$
are stable equilibria for the system \cref{wp}. In their linear stability analysis, Mori, Jilkine, and Edelstein-Keshet~\cite{Mori2011} observed that this stability condition can be reduced to the following:
\begin{equation}
    f_u(u_\pm(v),v) - f_v(u_\pm(v),v) < 0. \label{homogeneous-stability-condition}
\end{equation}

    \item[Velocity sign condition.] For the function
\begin{equation}
    I(v) = \int_{u_-(v)}^{u_+(v)}f(u,v)du, \label{I}
\end{equation}
there exists a unique $v_c \in (v_\mathrm{min}, v_\mathrm{max})$ such that $I(v_c) = 0$, with $I(v)$ being positive for $v > v_c$, negative for $v < v_c$. In addition, it holds that the derivative of $I$ at $v = v_c$ is strictly positive.
\end{description}
The velocity sign condition will be used in Section 3.
These conditions are indispensable kinetic requirements for the wave-pinning phenomenon.
Another key component for the wave-pinning model is the conservation of the total mass:
\begin{equation}
    \int_{\Omega}(u + v)dx = M, \label{mass-conservation}
\end{equation}
where
\[
M = \int_{\Omega}(u(x,0) + v(x,0))dx.
\]


Let $(u^\eps, v^\eps)$ be a solution to \cref{wp}, and suppose that $v^\eps \in [v_\mathrm{min}, v_\mathrm{max}]$.
On the fast timescale = $\tau = t/\eps$, the leading-order dynamics $\partial_\tau u^\eps = f(u^\eps,v^\eps)$ drive $u^\eps$ toward the stable states $u_\pm(v^\eps)$, leading to the rapid formation of a sharp interface $\Gamma^\eps(t) := \{x\in\Omega\mid u^\eps(x,t) = u_m(v^\eps)\}$ separating the regions $\{u^\eps \approx u_\pm(v^\eps)\}$.
In what follows, we will derive the equation of motion for the interface $\Gamma^\eps(t)$ that develops
on multiple timescales, based on matched asymptotic expansions.
Specifically, we employ the additional asymptotic expansion of the signed distance function, as described by Nakamura et al. \cite{Nakamura1999}, to derive higher-order approximations for the equation of motion of interfaces.
For fundamental concepts and procedures of the method of the matched asymptotic expansions, see, e.g., \cite{Caginalp1988,Ei1997,Nakamura1999,Rubinstein1989} and references therein.

\subsection{Matched asymptotic expansions}
We assume that the interface $\Gamma^\eps(t)$ is either a smooth closed curve or a smooth curve whose endpoints lie on the boundary $\partial\Omega$ for every $t \geq 0$.
The domain $\Omega$ is partitioned into two regions, $\Omega_-^\eps(t) := \{x \in \Omega\,|\,u^\eps < u_m(v^\eps)\}$ and
$\Omega_+^\eps(t) := \{x \in \Omega\,|\,u^\eps > u_m(v^\eps)\}$. Here, we define the signed distance function
to $\Gamma^\eps(t)$ as follows:
\begin{equation}
    d^\eps(x,t) = \begin{cases}
        \mathrm{dist}(x, \Gamma^\eps(t)), & x \in \Omega_-^\eps(t), \\
        0, & x \in \Gamma^\eps(t), \\
        -\mathrm{dist}(x, \Gamma^\eps(t)), & x \in \Omega_+^\eps(t),
    \end{cases}
\end{equation}
where $\mathrm{dist}(x, \Gamma^\eps(t))$ is the distance from $x \in \Omega$ to $\Gamma^\eps(t)$ in $\mathbb{R}^2$.
Then, we can rewrite
\begin{align*}
    &\Gamma^\eps(t) = \{x \in \Omega\,|\,d^\eps(x,t) = 0\}, & &\Omega_\pm^\eps(t) = \{x \in \Omega\,|\,\pm d^\eps(x,t) < 0\},
\intertext{and define}
    &\Gamma^\eps = \bigcup_{t \geq 0}(\Gamma^\eps(t) \times \{t\}), & &\Omega_\pm^\eps = \bigcup_{t \geq 0}(\Omega_\pm^\eps(t) \times \{t\}).
\end{align*}

We assume the following expansions for $u$ and $v$ in the outer region (i.e., away from the interface $\Gamma^\eps$):
\begin{subequations}
    \label{outer}
    \begin{align}
        u^\eps(x,t) &= u_0(x,t) + \eps u_1(x,t) + \eps^2 u_2(x,t) + \cdots, \label{outer-u} \\
        v^\eps(x,t) &= v_0(x,t) + \eps v_1(x,t) + \eps^2 v_2(x,t) + \cdots, \label{outer-v}
    \end{align}
\end{subequations}
where the functions $u_k$ and $v_k$ are smooth, bounded and independent of $\eps$ for $k = 0,1,2,\dots$.
These functions satisfy the no-flux boundary condition \cref{wp-bc}.
To be consistent with the total mass $M$ in \cref{mass-conservation}, we set
\begin{align}
    \int_{\Omega}(u_0 + v_0)dx &= M, \label{mass-conservation-outer0}\\
    \int_{\Omega}(u_k + v_k)dx &= 0, \quad k = 1,2,\dots. \label{mass-conservation-outer}
\end{align}

Let $z(x,t)$ be the stretched variable, defined by $z(x,t) = d^\eps(x,t) / \eps$.
Also, let $\gamma(s,t)$ be a parametrization of the interface $\Gamma^\eps(t)$,
and $n(s,t)$ be the unit normal vector to $\Gamma^\eps(t)$, pointing from $\Omega_+^\eps(t)$ to $\Omega_-^\eps(t)$.
Here, $s \in [0,L]$ is an arc-length parameter along $\Gamma^\eps(t)$, where $L \Add{= L(t)}$ is the length of $\Gamma^\eps(t)$.
Then, for any point $x$ in the neighborhood of $\Gamma^\eps(t)$, we have the unique expression
\[
x = \gamma(s(x,t),t) + \eps z(x,t)n(s(x,t),t).
\]
This allows us to introduce the local coordinate system $(s,z)$ near the interface $\Gamma^\eps(t)$.
Hence, we assume the following expansions for $u$ and $v$ in the inner region (i.e., near the interface $\Gamma^\eps$):
\begin{subequations}
    \label{inner}
    \begin{align}
        u^\eps(x,t) &= U^\eps(s,z,t) = U_0(s,z,t) + \eps U_1(s,z,t) + \eps^2 U_2(s,z,t) + \cdots, \label{inner-u} \\
        v^\eps(x,t) &= V^\eps(s,z,t) = V_0(s,z,t) + \eps V_1(s,z,t) + \eps^2 V_2(s,z,t) + \cdots. \label{inner-v}
    \end{align}
\end{subequations}
where the functions $U_k$ and $V_k$ are smooth and independent of $\eps$ for $k = 0,1,2,\dots$.

By matching the outer expansions \cref{outer-u}, \cref{outer-v} and the inner expansions \cref{inner-u}, \cref{inner-v},
we require that the following conditions hold for $s \in [0,L]$ and $t > 0$.
At the leading order, the large $z$ behavior of the inner solution must match the outer solution at the interface:
\begin{subequations}
    \begin{gather}
        U_0(s,\pm\infty,t) = u_0(\gamma(s,t) \pm 0, t) + o(1) \quad (z \to \pm\infty), \label{matching0-u}\\
        V_0(s,\pm\infty,t) = v_0(\gamma(s,t) \pm 0, t) + o(1) \quad (z \to \pm\infty). \label{matching0-v}
    \end{gather}
\end{subequations}
Here, $g(\gamma(s,t) \pm 0)$ denotes the limiting value of an outer function $g$ as $\Gamma^\eps(t)$ is approached from the side where $\pm d(x,t) > 0$, i.e., from $\Omega_\mp^\eps(t)$.
And at the first order, the large $z$ behavior of the inner solution must be linear in $z$:
\begin{subequations}
    \begin{gather}
        U_1(s,z,t) = u_1(\gamma(s,t) \pm 0, t) + z\partial_n u_0(\gamma(s,t) \pm 0, t) + o(1) \ (z \to \pm\infty), \label{matching1-u} \\
        V_1(s,z,t) = v_1(\gamma(s,t) \pm 0, t) + z\partial_n v_0(\gamma(s,t) \pm 0, t) + o(1) \ (z \to \pm\infty). \label{matching1-v}
    \end{gather}
\end{subequations}
Furthermore, at the second order, we have
\begin{subequations}
    \begin{gather}
        \begin{aligned}
        U_2(s,z,t) &= u_2(\gamma(s,t) \pm 0, t) + z\,\partial_n u_1(\gamma(s,t) \pm 0,t) \\
            &\quad + \frac{z^2}{2}\partial_n^2 u_0(\gamma(s,t) \pm 0, t) + o(1) \quad (z \to \pm\infty),
        \end{aligned} \label{matching2-u}\\
        \begin{aligned}
        V_2(s,z,t) &= v_2(\gamma(s,t) \pm 0, t) + z\,\partial_n v_1(\gamma(s,t) \pm 0,t) \\
            &\quad + \frac{z^2}{2}\partial_n^2 v_0(\gamma(s,t) \pm 0, t) + o(1) \quad (z \to \pm\infty).
        \end{aligned} \label{matching2-v}
    \end{gather}
\end{subequations}
For a derivation of these relations, see the appendix of \cite{Caginalp1988}.

We assume the following expansion for $d^\eps$:
\begin{equation}
    d^\eps(x,t) = d_0(x,t) + \eps d_1(x,t) + \eps^2 d_2(x,t) + \cdots,
\end{equation}
where each function $d_k$ is smooth, bounded and independent of $\eps$ for $k = 0,1,2,\dots$.

\subsection{Leading-order outer and inner solutions}
\paragraph{Outer solutions}
Substituting \cref{outer} into \cref{wp} and collecting the terms of order $\eps^0$, we obtain
\begin{subequations}
    \begin{align}
        f(u_0,v_0) &= 0, \label{outer-eps0-u}\\
        D\Delta v_0 - f(u_0,v_0) &= 0. \label{outer-eps0-v}
    \end{align}
\end{subequations}
From \cref{outer-eps0-u} and bistability of $f(\cdot,v_0)$, it holds that $u_0 = u_\pm(v_0)$ or $u_0 = u_m(v_0)$. Now we let
\begin{equation}
    u_0(x,t) = \begin{cases}
        u_-(v_0), & x \in \Omega_-^\eps(t), \\
        u_+(v_0), & x \in \Omega_+^\eps(t),
    \end{cases} \label{u0}
\end{equation}
at time $t > 0$.
Adding \cref{outer-eps0-u} and \cref{outer-eps0-v} yields $D\Delta v_0 = 0$ in $\Omega_-^\eps(t)$ and $\Omega_+^\eps(t)$,
which, together with \cref{wp-bc}, implies that $v_0$ is spatially uniform in $\Omega_-^\eps(t)$ and $\Omega_+^\eps(t)$. We denote
\begin{equation}
    v_0(x,t) = \begin{cases}
        v_<(t), & x \in \Omega_-^\eps(t), \\
        v_>(t), & x \in \Omega_+^\eps(t).
    \end{cases}
\end{equation}

\paragraph{Inner solutions}
Note that $|\nabla d^\eps| = 1$, $\nabla s \cdot \nabla d^\eps = 0$, $|\nabla s|^2 = 1 + O(\eps)$, $\Delta s = O(\eps)$,
and $\Delta d^\eps = \kappa - \eps\kappa^2 z + O(\varepsilon^2)$, where $\kappa = \kappa(s,t)$ is the curvature of $\Gamma^\eps(t)$.
\Add{Here, we adopt the sign convention that the curvature \( \kappa \) is positive if the region \( \Omega_+^\eps(t) \) is convex, with the interface parametrized such that \( \Omega_+^\eps(t) \) lies to its left.}
Using these relations and the chain rule, we obtain the following relations for \cref{inner-u} (and likewise for \cref{inner-v}) hold:
\begin{align*}
    \partial_t u^\eps &= \frac{1}{\eps}\,\partial_t d^\eps\Add{\cdot}\,\partial_z U^\eps + \partial_t s\,\Add{\cdot}\,\partial_s U^\eps + \partial_t U^\eps, \\
    \Delta u^\eps &= \frac{1}{\eps^2}\,\partial_z^2 U^\eps + \frac{1}{\eps}\kappa\,\partial_z U^\eps + \partial_s^2 U^\eps - \kappa^2 z\,\partial_z U^\eps + O(\eps).
\end{align*}
Thus, substituting \cref{inner-v} into \cref{wp-v} and collecting the terms of order $\eps^{-2}$, we therefore obtain
\begin{equation}
    D\partial_z^2 V_0 = 0, \label{inner-eps^(-2)-v}\\
\end{equation}
From \cref{inner-eps^(-2)-v}, $V_0$ is represented as
\[
V_0 = a(\Add{s,}t)z + b(\Add{s,}t),
\]
where $a$ and $b$ are real functions of \Add{\( s \) and}\Erase{time} $t$.
However, the matching condition \cref{matching0-v} implies that $V_0$ remains bounded as $z \to \pm\infty$,
leading to $a(\Add{s,}t) = 0$ and $b(\Add{s,}t) = v_<(t) = v_>(t)$.
Therefore, $V_0$ reduces to \Erase{the}\Add{a} spatially uniform function $v_0(t)$:
\begin{equation}
    V_0(\Add{s,z,}t) = v_0(t).
\end{equation}
\Add{The above discussion also implies that the leading-order outer solution for \( v \) is continuous across \( \Gamma^\eps(t) \) and takes the constant value \( v_0(t) \) in all of \( \Omega \).}
We derive the governing equation for $v_0$.
Integrating $u_0 + v_0$ over $\Omega$ with \cref{u0} yields:
\[
\int_{\Omega}(u_0 + v_0)dx = u_+(v_0)|\Omega_+^\eps(t)| + u_-(v_0)|\Omega_-^\eps(t)| + v_0|\Omega|,
\]
where $|\Omega|$ represents the area of $\Omega$ and similarly for the subdomains.
Thus, the mass conservation law \cref{mass-conservation-outer0} reduces to
\begin{equation}
    u_+(v_0)|\Omega_+^\eps(t)| + u_-(v_0)|\Omega_-^\eps(t)| + v_0|\Omega| = M, \label{reduced-mass-conservation}
\end{equation}
which determines $v_0$.

Next, we consider the inner solutions for $u$. Substituting \cref{inner-u} into \cref{wp-u} and collecting the terms of order $\eps^0$, we obtain
\begin{equation}
    \partial_z^2 U_0 - \partial_t d_0\,\Add{\cdot}\,\partial_z U_0 + f(U_0,V_0) = 0. \label{inner-eps^0-u}
\end{equation}
From the matching condition \cref{matching0-u}, $U_0$ is characterized as a heteroclinic solution from $u_+(V_0)$ toward $u_-(V_0)$.
Multiplying \cref{inner-eps^0-u} by $\partial_z U_0$ and integrating over $z \in \mathbb{R}$, we obtain
\begin{equation}
    -\partial_t d_0 = c(V_0) := \frac{\int_{u_-(V_0)}^{u_+(V_0)}f(u,V_0)du}{\int_{\mathbb{R}}(\partial_z U_0)^2dz}. \label{normal-velocity0}
\end{equation}

Note that $U_0$ depends only on $z$ and $V_0$. Therefore, we write $U_0 = U_0(z,V_0)$.

\subsection{First-order outer and inner solutions}
\paragraph{Outer solutions}
Substituting \cref{outer-u} into \cref{wp-u} and collecting the terms of order $\eps^1$, we obtain
\begin{subequations}
    \begin{align}
        \partial_t u_0 &= f_u(u_0,v_0)u_1 + f_v(u_0,v_0)v_1, \label{outer-eps1-u}\\
        \partial_t v_0 &= D\Delta v_1 - f_u(u_0,v_0)u_1 - f_v(u_0,v_0)v_1. \label{outer-eps1-v}
    \end{align}
\end{subequations}
Adding \cref{outer-eps1-u} and \cref{outer-eps1-v}, we obtain
\begin{equation}
    D\Delta v_1 = \partial_t (u_0 + v_0). \label{v1-poisson}
\end{equation}

We now derive the solvability condition for the Poisson equation \cref{v1-poisson} by integrating it over the entire domain $\Omega$.
First, we evaluate the left-hand side. Applying the divergence theorem and using the no-flux boundary condition for $v_1$ on $\partial\Omega$,
the integral of the Laplacian term is reduced to an integral over the interface $\Gamma^\eps(t)$:
\begin{equation}
    \int_{\Omega}\Delta v_1 dx = \int_{\Gamma^\eps(t)}\llbracket \partial_n v_1 \rrbracket ds.
\end{equation}
Here, $\llbracket g \rrbracket$ denotes the jump of a quantity $g = g(x,t)$ across the interface $\Gamma^\eps(t)$, defined as $\jump{g} := g(\gamma(s,t) - 0, t) - g(\gamma(s,t) + 0, t)$
for $s \in [0,L]$ and $t > 0$.

Next, we evaluate the right-hand side. The total mass of the leading order terms is conserved as in \cref{mass-conservation-outer0}.
Applying the Reynolds transport theorem relates this conservation law to the source term:
\[
0 = \frac{d}{dt}\int_{\Omega}(u_0 + v_0)dx = \int_{\Omega}\partial_t(u_0 + v_0)dx + \jump{u_0}\int_{\Gamma^\eps(t)}V_n\,ds,
\]
where $V_n = V_n(s,t)$ is the normal velocity of the interface $\Gamma^\eps(t)$.
From this equation, we have
\begin{equation}
    \int_{\Omega}\partial_t(u_0 + v_0)dx = -\jump{u_0}\int_{\Gamma^\eps(t)}V_n\,ds. \label{reynolds}
\end{equation}
Therefore, as a solvability condition for \cref{v1-poisson}, we conclude that
\begin{equation}
    \int_{\Gamma^\eps(t)}\big(D\jump{\partial_n v_1} + \jump{u_0}V_n\big)ds = 0. \label{solvability-condition-poisson}
\end{equation}

\color{black}We also need to rewrite the right-hand side of the Poisson equation \cref{v1-poisson} with another form, which is $\partial_t(u_0 + v_0)$ now.
We rewrite \cref{reynolds} as
\begin{equation}
    \Lambda(v_0)\Erase{\partial_t v_0}\Add{\dot{v}_0} = -\jump{u_0}\dot{A}, \label{reynolds-simple}
\end{equation}
where $\dot{A} := \int_{\Gamma^\eps(t)}V_n\,ds$ is the time derivative of the area of the subdomain $\Omega_+^\eps(t)$,
\begin{equation}
    \Lambda(v_0) := u_+'(v_0)|\Omega_+^\eps(t)| + u_-'(v_0)|\Omega_-^\eps(t)| + |\Omega| \label{Lambda}
\end{equation}
and the prime ($'$) denotes the derivative with respect to $v_0$.
\Add{Here and hereinafter, the dot ( \(\dot{}\) ) denotes the derivative with respect to \( t \) (i.e., \( \dot{v}_0 = dv_0 / dt \)).}

We claim that $\Lambda(v_0)$ is positive.
To show this, we differentiate $f(u_\pm(v),v) = 0$ with respect to $v$ to obtain $f_u(u_\pm(v),v)u_\pm'(v) + f_v(u_\pm(v),v) = 0$.
The stability condition \cref{homogeneous-stability-condition} ensures that $f_v(u_\pm(v),v) > f_u(u_\pm(v),v)$, which leads to the inequality
\[
(1 + u_\pm'(v))f_u(u_\pm(v),v) < 0.
\]
For this product to be negative under \cref{bistability-condition}, we must conclude that $1 + u_\pm'(v) > 0$, which confirms the positivity of $\Lambda(v_0)$.
Combined with \cref{reynolds-simple}, this implies that \Erase{$\partial_t v_0$}\Add{\(\dot{v}_0\)} and $\dot{A}$ have opposite signs,
that is, $v_0$ decreases as the subdomain $\Omega_+^\eps(t)$ expands.

We also calculate $\partial_t u_0 = u_\pm'(v_0)\Erase{\partial_t v_0}\Add{\dot{v}_0}$ in $\Omega_\pm^\eps$.
Thus, $\partial_t(u_0 + v_0) = (1 + u_\pm'(v_0))\Erase{\partial_t v_0}\Add{\dot{v}_0}$ holds and the Poisson equation \cref{v1-poisson} can be rewritten into
\begin{equation}
    -\Delta v_1 = \frac{(1 + u_\pm'(v_0))\jump{u_0}}{D\Lambda(v_0)}\dot{A} \qquad \text{in } \Omega_\pm^\eps \label{v1-poisson-simple}
\end{equation}
by using \cref{reynolds-simple}.

\paragraph{Inner solutions}
Substituting \cref{inner-v} into \cref{wp-v} and collecting the terms of order $\eps^{-1}$, we obtain
\begin{gather}
    D(\partial_z^2 V_1 + \kappa\,\partial_z V_0) = 0. \label{inner-eps^(-1)-v}
\end{gather}
Since $\partial_z V_0 = 0$, we obtain $D\partial_z^2 V_1 = 0$. By an argument analogous to the one for $V_0$,
we conclude that $v_1$ is continuous on the interface $\Gamma^\eps(t)$ and $V_1$ is constant in $z$:
\begin{equation}
    V_1(s,z,t) = v_1(\gamma(s,t),t), \qquad s \in [0,L],\ t > 0.
\end{equation}

Substituting \cref{inner-u} into \cref{wp-u} and collecting the terms of order $\eps^1$, we obtain
\begin{equation}
    \begin{aligned}
        &\partial_z^2 U_1 - \partial_t d_0\,\Add{\cdot}\,\partial_z U_1 + f_u(U_0,V_0)U_1 \\
        =\ &\partial_t s\,\Add{\cdot}\,\partial_s U_0 \,+\partial_tU_0 + \partial_t d_1\Add{\cdot}\,\partial_z U_0 - \kappa\,\partial_z U_0 - f_v(U_0,V_0)V_1.
    \end{aligned} \label{inner-eps^1-u}
\end{equation}
Let $\mathcal{L}$ be an operator on $L^2(\mathbb{R})$ defined by $\mathcal{L} := \partial_z^2 - (\partial_t d_0)\partial_z + f_u(U_0,V_0)$,
and $h \in L^2(\mathbb{R})$ be the right-hand side of \cref{inner-eps^1-u}. Then, we can rewrite \cref{inner-eps^1-u} as
$\mathcal{L}U_1 = h$, i.e., $h$ is in the range of $\mathcal{L}$. Differentiating \cref{inner-eps^0-u} with respect to $z$, we also see that $\mathcal{L}(\partial_z U_0) = 0$.
Thus, for the adjoint operator of $\mathcal{L}$ denoted by $\mathcal{L}^* = \partial_z^2 + (\partial_t d_0)\partial_z + f_u(U_0,V_0)$,
$p^*(z) = \Erase{\partial_z U_0(-z)}\Add{e^{-cz}\partial_z U_0(z)}$ \Add{with \( c = \partial_t d_0 \) }satisfies $\mathcal{L}^* p^* = 0$, i.e., $p^*$ is in the kernel of $\mathcal{L}^*$.
Since the range of $\mathcal{L}$ is orthogonal to the kernel of its adjoint $\mathcal{L}^*$, we have
\[
\int_{\mathbb{R}}\left(\partial_t U_0 + (\partial_t d_1 - \kappa)\partial_z U_0 - f_v(U_0,V_0)V_1\right)p^*dz = 0,
\]
where we used the fact that $\partial_s U_0 = 0$. Therefore, we conclude that
\begin{equation}
    \begin{aligned}
        -\partial_t d_1 &= -\kappa + \frac{\int_{\mathbb{R}}\left(\partial_t U_0 - f_v(U_0,V_0)V_1\right)\Add{e^{-cz}\partial_z U_0}\Erase{\partial_z U_0(-z)}dz}{\int_{\mathbb{R}}\Add{e^{-cz}}\partial_z U_0^{\Add{2}}\Erase{(z)\partial_z U_0(-z)}dz} \\
        &= -\kappa + \alpha(V_0) - \beta(V_0)V_1,
    \end{aligned} \label{normal-velocity1}
\end{equation}
where
\begin{equation}
    \alpha(V_0) = \frac{\int_{\mathbb{R}}\Add{e^{-cz}}\partial_t U_0 \Add{\partial_z U_0}\Erase{\partial_z U_0(-z)}dz}{\int_{\mathbb{R}}\Add{e^{-cz}}\partial_z U_0^{\Add{2}}\Erase{(z)\partial_z U_0(-z)}dz}, \quad \beta(V_0) = \frac{\int_{\mathbb{R}}\Add{e^{-cz}}f_v(U_0,V_0)\Add{\partial_z U_0}\Erase{\partial_z U_0(-z)}dz}{\int_{\mathbb{R}}\Add{e^{-cz}}\partial_z U_0^{\Add{2}}\Erase{(z)\partial_z U_0(-z)}dz}. \label{B(v_0)}
\end{equation}

\subsection{The interface equation}
Using the relation $V_n = -\partial_t d^\eps = -\partial_t d_0 - \eps\partial_t d_1 + O(\eps^2)$,
\cref{normal-velocity0} and \cref{normal-velocity1} yield the following approximate interface equation for $\Gamma^\eps$:
\begin{equation}
    V_n = c(v_0) - \eps\left(\kappa - \alpha(v_0) + \beta(v_0)v_1\right) \qquad \text{on } \Gamma^\eps. \label{interface-equation}
\end{equation}
We note that the Poisson equation \cref{v1-poisson-simple} for $v_1$ has a unique solution subject to the no-flux boundary condition \cref{wp-bc}, the mass conservation law \cref{mass-conservation-outer},
and continuity across the interface $\Gamma^\eps(t)$.
Here, the term $u_1$ appearing in the conservation law \cref{mass-conservation-outer} is explicitly obtained from \cref{outer-eps1-u}.

If the interface touches the boundary $\partial\Omega$ at its endpoints, the no-flux boundary condition \cref{wp-bc} ensures that the interface meets the boundary $\partial\Omega$ orthogonally \cite{Rubinstein1989}.

Lastly, we shall simplify the interface equation \cref{interface-equation}.
To this end, we begin by deriving the pointwise jump relation between $u_0$ and $\partial_n v_1$ from the matching of outer and inner solutions.
Integrating \cref{inner-eps^0-u} over $\mathbb{R}$ and using \cref{matching0-u}, we obtain
\begin{equation}
    \int_{\mathbb{R}}f(U_0,v_0)dz = c(v_0)\jump{u_0}. \label{jump-u0}
\end{equation}
Next, substituting \cref{inner-v} into \cref{wp-v} and collecting the term of order $\eps^0$, we obtain
\begin{equation}
    f(U_0,v_0) = D\,\partial_z^2 V_2.
\end{equation}
Integrating this equation over $\mathbb{R}$ and using the second-order matching condition \cref{matching2-v}, we obtain
\begin{equation}
    \int_{\mathbb{R}}f(U_0,v_0)dz = -D\jump{\partial_n v_1}. \label{jump-v1}
\end{equation}
Therefore, \cref{jump-u0} and \cref{jump-v1} lead to the following relation:
\begin{equation}
    D\jump{\partial_n v_1} + c(v_0)\jump{u_0} = 0. \label{jump-relation}
\end{equation}
We now compare this jump relation with the solvability condition \cref{solvability-condition-poisson} for $v_1$,
discarding terms of order $O(\eps^2)$ of $V_n$, yielding
\begin{equation}
    \jump{u_0}\int_{\Gamma^\eps(t)}\left(-\kappa + \alpha(v_0) - \beta(v_0)v_1\right)ds = 0. \label{alpha-beta}
\end{equation}
Since $\alpha(v_0)$ and $\beta(v_0)$ are spatially constants, \cref{alpha-beta} gives the expression for $\alpha(v_0)$:
\begin{equation}
    \alpha(v_0) = \average{\kappa} + \beta(v_0)\average{v_1}, \label{A(v_0)}
\end{equation}
where $\average{\kappa}$ is the average of $\kappa$ over the interface $\Gamma^\eps(t)$ and likewise for $\average{v_1}$.
Consequently, we arrive at the final form of the interface equation:
\begin{equation}
    V_n = c(v_0) - \eps\left[\kappa - \average{\kappa} + \beta(v_0)(v_1 - \average{v_1})\right]. \label{interface-equation2}
\end{equation}
Equation \cref{interface-equation2} encapsulates the interplay of forces acting on the interface. The term proportional to $\kappa$ represents the surface tension effect that drives the curvature flow, minimizing the interface length. The non-local average term $\average{\kappa}$ acts as a global Lagrange multiplier that enforces the area preservation, counteracting the shrinking effect of the curvature flow. The deviation term involving $v_1$ represents higher-order mass transport effects.
By virtue of this consequence of the area preservation on the slow timescale, we obtain that $\dot{A} = c(v_0)L$, allowing us to rewrite \cref{v1-poisson} as \cref{fbp-vp}.

This resulting equation reveals that the interface dynamics involve two distinct timescale effects:
a dominant $O(1)$ motion propagating at speed $c(v_0)$, and a slower $O(\eps)$ modulation driven by area-preserving curvature flow and higher-order terms.

We summarize our analysis in the following.
\begin{result}
    Assume the reaction term $f$ satisfies \cref{bistability-condition} and \cref{homogeneous-stability-condition}.
    Let $(u^\eps, v^\eps)$ be a solution to the wave-pinning model \cref{wp} consistent with the matched asymptotic expansions defined in \cref{outer} and \cref{inner}.
    Then, the dynamics of the interface $\Gamma^\eps$ of the solution $u^\eps$ is approximated by the free boundary problem \cref{fbp} up to $O(\eps^2)$.
\end{result}

\section{Slow interface motion along the boundary}
In this section, we show that the interface, which intersects the boundary orthogonally and evolves according to \cref{interface-equation2},
moves toward a local maximum of the boundary curvature.
Our analysis employs a perturbation method for an interface enclosing a small area, following the approach of Alikakos et al. \cite{Alikakos2000}

Let $\Omega \subset \mathbb{R}^2$ be a bounded domain with smooth boundary $\partial\Omega$ of class $C^4$.
We consider a family of smooth simple curves $\{\Gamma(t)\}_{t \geq 0}$, which evolve according to the interface equation \cref{interface-equation2}.
Each curve $\Gamma(t)$ intersects the boundary $\partial\Omega$ orthogonally at its endpoints, and separates the domain into $\Omega_+(t)$ and $\Omega_-(t)$.
We assume that the evolution of the interface is in a quasi-static regime. More precisely, we consider the case where the interface evolves slowly with a speed of $O(\eps)$, or equivalently, the magnitude of the leading-order velocity $c(v_0)$ is $O(\eps)$.
In such a slow regime, the interface tends to remain close to a stable equilibrium state.
Indeed, numerical simulations (Fig. \ref{fig:long-time-behavior}) suggest that the interface evolves while maintaining a shape sufficiently close to a circular arc.
Based on this observation, we assume that $\Gamma(t)$ is close to a circular arc.
To localize our analysis, we also assume that the size of the nearly circular interface is sufficiently small; specifically, we consider the case where the interface length $L \Add{= L(t)}$ is of order $\eps$.
Under these assumptions, the free boundary problem \cref{fbp} reduces to the following simplified form, independent of $v_1$.
\begin{subequations}
    \label{reduced-interface-system}
    \begin{align}
        V_n = c(v_0) - \eps(\kappa - \average{\kappa}) + O(\eps^3) && \text{on }\Gamma(t),\ t > 0, \label{slow-interface-equation}\\
        u_+(v_0)|\Omega_+(t)| + u_-(v_0)|\Omega_-(t)| + v_0|\Omega| = M && t \geq 0. \label{slow-mass-conservation}
    \end{align}
\end{subequations}
\Erase{Here, $v_1$ is assumed to be small and is absorbed into the $O(\eps^3)$ term.
While a rigorous derivation of the $L^\infty$ bounds for $v_1$ requires a detailed analysis of the transmission problem with jump condition \cref{jump-relation}, the scaling behavior can be understood as follows. Recall that $v_1$ satisfies the Poisson equation \cref{v1-poisson-simple}, where the source term is proportional to $c(v_0)L$. Given the assumptions that $c(v_0)$ and $L$ are of order $\eps$, the product $c(v_0)L$ scales as $O(\eps^2)$. Thus, standard elliptic regularity estimates suggest that the magnitude of $v_1$ is controlled by this source term, justifying the ansatz $v_1=O(\eps^2)$.}
\Erase{For simplicity, we rename $T$ back to $t$ and hereafter}
\Add{Here the term involving $v_{1}$ is absorbed into the $O(\eps^{3})$ remainder.
We emphasize that the transmission problem \cref{fbp-vp,fbp-jump,fbp-nf} determines $v_{1}$ only up to an additive constant; this constant is fixed by the mass constraint \cref{fbp-mc1}.
Since the interface equation contains $v_{1}$ only through $v_{1}-\langle v_{1}\rangle$, the relevant estimate for the present reduction is
$$v_{1}-\langle v_{1}\rangle = O(\eps^{2}) \quad \text{on } \Gamma.$$
Let us explain the formal scaling behind this estimate. In the regime considered here, $c(v_{0})=O(\eps)$ and $L=O(\eps)$. Hence the bulk source term in \cref{fbp-vp} is of order $O(\eps^{2})$, while the jump condition \cref{fbp-jump} gives a normal derivative jump of order $O(\eps)$. After rescaling a neighborhood of the small interface by $x=X_{0}(s_{0})+\eps y$, the corresponding scaled problem has a bulk source of order $O(\eps^{4})$ and a conormal derivative jump of order $O(\eps^{2})$. Standard Lipschitz and gradient estimates for elliptic transmission problems with piecewise regular data, such as those in \cite{dong2019gradient,li2000gradient}, therefore support the formal scaling
$$\|\nabla v_{1}\|_{L^\infty(\Gamma')} = O(\eps)$$
on compact portions $\Gamma'$ of the interface away from the boundary contact points. Consequently, 
\[\|v_{1}-\langle v_{1}\rangle\|_{L^{\infty}(\Gamma)}\le CL\|\partial_{s}v_{1}\|_{L^{\infty}(\Gamma)}=O(\eps^{2}).\]
Thus $\eps\beta(v_{0})(v_{1}-\langle v_{1}\rangle)=O(\eps^{3})$, which justifies its absorption into the remainder in \cref{slow-interface-equation} at the formal level.}

\Add{We note that obtaining an estimate uniform up to the points where $\Gamma$ meets $\partial\Omega$ is more delicate. Boundary-intersecting transmission problems generally require additional compatibility analysis: see, for example, \cite{druet2013global}. A rigorous justification of this estimate in the present $\eps$-dependent geometry is left for future work.}

In what follows, we introduce a conformal mapping that flattens the curved boundary, thereby transferring our analysis to the complex upper half-plane. In this setting, we treat the interface as a perturbation of a semi-circle and compute the relevant geometric quantities, such as the normal velocity and mean curvature. Finally, applying these to the reduced free boundary problem, we demonstrate that the interface moves along the boundary driven by the gradient of the boundary curvature, and we establish the stability of the near-circular configuration.

\subsection{Conformal mapping}
\Add{Hereafter, we identify $\mathbb{R}^2$ with $\mathbb{C}$.}
Analyzing interface motion near a general curved boundary is analytically challenging due to the complexity of the Neumann boundary condition. To overcome this, we employ a conformal mapping technique. This allows us to ``flatten'' the boundary geometry, transforming the complex domain interactions into a tractable perturbation problem on the upper half-plane.

Let $X^0(s),\,s\in\mathbb{R}/L_{\Add{\partial\Omega}}\mathbb{Z}$, be an arc-length parametrization of $\partial\Omega$, where $L_{\Add{\partial\Omega}}$ is the length of $\partial\Omega$.
Since the boundary $\partial\Omega$ is of class $C^4$, $X^0$ is also in $C^4(\mathbb{R}/L_{\Add{\partial\Omega}}\mathbb{Z})$.
Let $z_0 \in \Omega$ be an arbitrary point. We denote the upper half-plane of $\mathbb{C}$ as $\mathbb{H}$.
Then, there exists a unique conformal diffeomorphism $\Psi \colon \overline{\mathbb{H}} \times \mathbb{R}/L_{\Add{\partial\Omega}}\mathbb{Z} \to \overline{\Omega}$
satisfying the following properties for every $s \in \mathbb{R}/L_{\Add{\partial\Omega}}\mathbb{Z}$:
\begin{itemize}
    \item $\Psi(\cdot,s)$ maps the real line onto $\partial\Omega$ and a point on the positive imaginary axis to the point $z_0$.
    \item $\Psi(0,s) = X^0(s)$, $\Psi_\zeta(0,s) = \tau(s) := X_s^0(s)$.
\end{itemize}
This conformal mapping is also of class $C^4$, and its existence and uniqueness are guaranteed by the Riemann mapping theorem.
The overall approach for this mapping, including the derivation of its properties, is detailed in Alikakos et al. \cite{Alikakos2000}.

For any fixed $s_0 \in \mathbb{R}/L_{\Add{\partial\Omega}}\mathbb{Z}$, let $g(s,s_0)$ be the real-valued smooth function satisfying
\begin{equation}
    X^0(s + s_0) = \Psi(g(s,s_0), s_0), \qquad s \in (-L_{\Add{\partial\Omega}}/4,L_{\Add{\partial\Omega}}/4). \label{X^0-Psi}
\end{equation}
From the definition of $\Psi$, we have $g(0,s_0) = 0$. The Taylor expansion of $\Psi(g(s,s_0),s_0)$ around $s = 0$ yields
\begin{equation}
    \begin{aligned}
        \Psi(g(s,s_0), s_0) =& X^0(s_0) + \Psi_\zeta g_s s + \frac{1}{2}\left(\Psi_{\zeta\zeta} g_s^2 + \Psi_\zeta g_{ss}\right)s^2 \\
        &+ \frac{1}{6}\left(\Psi_{\zeta\zeta\zeta} g_s^3 + 3\Psi_{\zeta\zeta} g_{ss} g_s + \Psi_\zeta g_{sss}\right)s^3 + O(s^4),
    \end{aligned} \label{taylor-Psi}
\end{equation}
where all the derivatives are evaluated at $\zeta = 0$ and at $s = 0$.
The Taylor expansion of $X^0(s + s_0)$ around $s = 0$ also yields
\begin{equation}
        X^0(s + s_0) = X^0(s_0) + \tau s + \frac{1}{2}iK\tau s^2 + \frac{1}{6}\left(iK_s - K^2\right)\tau s^3 + O(s^4), \label{taylor-X0}
\end{equation}
where we have used the relations $X_s^0 = \tau$, $X_{ss}^0 = \tau_s = -K\nu = iK\tau$, and $X_{sss}^0 = i(K\tau)_s = i(K_s\tau + iK^2\tau)$.
Here, $K$ is the curvature and $\nu$ is the unit outward normal vector of the boundary $\partial\Omega$. $i$ is the imaginary unit.
All terms on the right-hand side of \cref{taylor-X0} are evaluated at $s = 0$.

Comparing the coefficients of the expansions \cref{taylor-Psi} and \cref{taylor-X0} according to \cref{X^0-Psi}, we deduce that
$g_s(0,s_0) = 1$, and
\begin{align*}
    \Psi_{\zeta\zeta}(0,s_0) &= \left(iK(s_0) - g_{ss}(0,s_0)\right)\tau(s_0) =: a(s_0)\tau(s_0), \\
    \Psi_{\zeta\zeta\zeta}(0,s_0) &= \left(i(K_s(s_0) \!-\! 3K(s_0)g_{ss}(0,s_0)) + 3g_{ss}^2(0,s_0) \!-\! g_{sss}(0,s_0) \!-\! K^2(s_0)\right)\tau(s_0) \\
    &=: b(s_0)\tau(s_0).
\end{align*}

Since $\Psi$ is of class $C^4$ up to the real axis, for sufficiently small $\zeta$, we conclude that
\begin{equation}\label{eq:taylor-psi}
    \Psi(\zeta,s) = X^0(s) + \tau(s)\left(\zeta + \frac{1}{2}a(s)\zeta^2 + \frac{1}{6}b(s)\zeta^3 + O(\zeta^4)\right).
\end{equation}
Writing $a(s) = a_1(s) + ia_2(s)$ and $b(s) = b_1(s) + ib_2(s)$, where $a_k$ and $b_k$ are real-valued function for $k = 1,2$,
we see that
\begin{gather}\label{eq:psi-ab}
\begin{aligned}
    a_1(s) &= -g_{ss}(0,s), & a_2(s) &= K(s), \\
    b_1(s) &= 3g_{ss}^2(0,s) - g_{sss}(0,s) - K^2(s), & b_2(s) &=K_s(s) - 3K(s)g_{ss}(0,s).
\end{aligned}
\end{gather}

\subsection{Representation of the interface}
As mentioned at the beginning of this section, we consider the situation in which the interface $\Gamma(t)$ is sufficiently close to a circular arc with a perimeter of order $\eps$, which is small enough to allow the boundary $\partial\Omega$ to be locally approximated as flat on the scale of $\Gamma(t)$.
Specifically, we parametrize the interface $\Gamma(t)$ as follows:
\begin{equation}
    \Gamma(t) = \{\Psi(\zeta,s)\,|\,\zeta = \eps e^{\eps^2 R(\theta,t) + i\theta},\,s = s(t),\,0\leq\theta\leq\pi\}, \label{interface-representation}
\end{equation}
where the unknown $R = R(\theta,t)$ is a smooth function representing the perturbation,
and the unknown $s(t)$ is smooth and describes the center of the semicircle.
This parametrization effectively reduces the dynamics of $\Gamma(t)$ along the boundary to the dynamics of $s(t)$.

Since the interface touches the boundary orthogonally, we expand the perturbation $R(\theta,t)$ as a Fourier cosine series:
\begin{equation}
    R(\theta,t) = R_0(t) + \sum_{n \geq 2}R_n(t)\cos(n\theta). \label{R-expansion}
\end{equation}
\Add{Provided \( \eps > 0 \) is sufficiently small,} the $n = 1$ mode corresponds to a shift of the center, \Add{since the leading term of the barycenter of the droplet is proportional to the \( n=1 \) mode of \( R(\theta,t) \).}
\Add{Therefore, this mode is omitted from the expansion so that its effect is absorbed into $s(t)$.}

\subsection{Propagating wave speed \texorpdfstring{$c(v_0)$}{c(v0)}}
\label{condition for w}
Based on the parametrization for $\Gamma(t)$ as in \cref{interface-representation}, we write $v_0$ and $c(v_0)$ in terms of $R_0$.
As a first step, we calculate the area of the droplet enclosed by $\Gamma(t)$ and a segment of $\partial\Omega$.
We now let $\Omega_+(t)$ correspond to the droplet.
One can also choose $\Omega_-(t)$ as the droplet; the choice is of no consequence to the following analysis except for the signs of geometric quantities.
Let $D_+$ be the domain in the mapped coordinates, $D_+ := \{re^{i\theta}\,|\,0 \leq r \leq \eps e^{\eps^2 R(\theta,t)},\,0\leq\theta\leq\pi\}$.
The area is then:
\begin{equation}
\begin{aligned}\label{eq:area}
    |\Omega_+(t)| =& |\Psi(D_+,s)| = \int_{\Psi(D_+,s)}dw = \int_{\Erase{D^+}\Add{D_+}}|\Psi_\zeta|^2d\zeta \\
    =& \int_{D_+}|1 + a\zeta + b\zeta^2/2 + O(\zeta^3)|^2d\zeta \\
    =& \int_{0}^{\pi}\int_{0}^{\eps e^{\eps^2 R}}|1 + are^{i\theta} + br^2 e^{2i\theta}/2 + O(r^3)|^2rdrd\theta \\
    =& \frac{\pi}{2}\eps^2 - \frac{4}{3}\eps^3K + \pi\eps^4\left[R_0 + \frac{1}{4}(a_1^2 + a_2^2)\right] + O(\eps^5).
\end{aligned}
\end{equation}
Therefore, the area can be considered a function in terms of $R_0$.

Next, we denote the left-hand side of \cref{slow-mass-conservation} by $F(v_0, R_0)$.
The partial derivative of $F$ with respect to $v_0$ is given by $\Lambda(v_0)$ in \cref{Lambda}.
Since $\Lambda(v_0)$ is positive, the Implicit Function Theorem guarantees the existence of a unique function $v_0 = v_0(R_0)$
satisfying $F(v_0(R_0), R_0) = M$.
In addition, we have
\begin{align*}
    F_{R_0} &= u_+(v_0)\partial_{R_0}|\Omega_+(t)| + u_-(v_0)\partial_{R_0}(|\Omega| - |\Omega_+(t)|) \\
    &= (u_+(v_0) - u_-(v_0))(\pi\eps^4 + O(\eps^5)) > 0,
\end{align*}
leading to
\begin{equation}
    v_0'(R_0) = -\frac{F_{R_0}}{F_{v_0}} \Add{= -\frac{u_+(v_0(R_0)) - u_-(v_0(R_0))}{\Lambda(v_0(R_0))}(\pi\eps^4 + O(\eps^5))} < 0.
\end{equation}
Here, the prime ($'$) also denotes the derivative with respect to $R_0$. \Add{Note that \( v_0'(R_0) \) is of order \( \eps^4 \) since \( \Lambda(v_0) = (1 + u_-'(v_0))|\Omega| + O(\eps^2) \) holds from \cref{Lambda} and \cref{eq:area}.}

To be consistent with the assumption that \Erase{$c(v_0)$ is of order $\eps$}\Add{\(c(v_0) = O(\eps)\)}, we introduce a function $w_\eps(R_0)$ \Erase{of order $\eps^0$} such that
\begin{equation}
    v_0(R_0) = v_c + \eps w_\eps(R_0),
\end{equation}
where $v_c$ is a constant satisfying $c(v_c) = 0$ and $c'(v_c) > 0$ (see the velocity sign condition).
\Add{We deduce that \( w_\eps'(R_0) \) is of order \( \eps^3 \), since \( v_0'(R_0) = \varepsilon w_\eps'(R_0) \).}
\Add{Let \( R_0^0 \) be a value such that \( w_\eps(R_0^0) = 0 \), i.e., \( v_0(R_0^0) = v_c \), which leads to the representation \( w_\eps(R_0) = \int_{R_0^0}^{R_0}w_\eps'(r)dr \).}
\Add{Thus, it follows that \( w_\eps(R_0) \) is also of order \( \eps^3 \).}
Accordingly, we expand $c(v_0)$ into a Taylor series to obtain:
\begin{equation}
    c(v_0) = \eps c'(v_c)w_\eps(R_0) + \Erase{O(\eps^2 w_\eps^2)}\Add{O(\eps^{8})}. \label{traveling-wave-speed}
\end{equation}
It is clear that $w_\eps(R_0)$ also satisfies $w_\eps'(R_0) < 0$. \Add{From the above discussion, it also follows that \( \eps^{-3}w_\eps = O(1) \), which does not vanish at \( \eps = 0 \).}

\subsection{Geometric quantities}
We now calculate the normal velocity $V_n$ and curvature $\kappa$ of the interface $\Gamma(t)$, which is parametrized by \cref{interface-representation}.
We denote a point on the interface by $x(\theta,t) + iy(\theta,t) = z(\theta,t) = \Psi(\zeta,s(t))$, where $\zeta = \eps e^{\eps^2 R(\theta,t) + i\theta}$.
The unit normal vector is taken to point from the droplet region $\Omega_+(t)$ to $\Omega_-(t)$.

First, we calculate the normal velocity $V_n$. The normal velocity is given by
\begin{align*}
    V_n(\theta,t) &= \frac{x_t y_\theta - x_\theta y_t}{|z_\theta|} = -\frac{\Image(z_t \bar{z}_\theta)}{|z_\theta|} \\
    &= -\frac{\Image((\Psi_\zeta \zeta_t + \Psi_s \dot{s})\bar{\Psi}_\zeta\bar{\zeta}_\theta)}{|\Psi_\zeta \zeta_\theta|} = -\frac{\Image((|\Psi_\zeta|^2 \zeta_t + \Psi_s \bar{\Psi}_\zeta \dot{s})\bar{\zeta}_\theta)}{|\Psi_\zeta \zeta_\theta|}.
\end{align*}
\Erase{Here, the dot ( $\dot{}$ ) denotes the derivative with respect to time $t$.}
We see that
\begin{align*}
    \Psi_\zeta &= \tau + O(\zeta),& \Psi_s &= \tau + O(\zeta), \\
    \zeta_\theta &= (\eps^2 R_\theta + i)\zeta, & \Add{\zeta_t} &= \eps^2 \dot{R}\zeta.
\end{align*}
Substituting these into the formula for $V_n$ yields
\begin{equation}
    V_n(\theta,t) = \dot{s}\left(\cos\theta + \eps B_1^\eps\right) + \eps^3 \dot{R} \left(1 + \eps B_2^\eps\right) \label{normal-velocity}
\end{equation}
where $B_i^\eps$ is a function of order $\eps^0$ for $i = 1,2$.

Next, we calculate the curvature $\kappa$ of $\Gamma(t)$. The curvature is given by
\[
    \kappa(\theta,t) = \frac{y_{\theta\theta}x_\theta - x_{\theta\theta}y_\theta}{|z_\theta|^3} = \frac{\Image(z_{\theta\theta}\bar{z}_\theta)}{|z_\theta|^3} 
    = \frac{\Image(\zeta_{\theta\theta}\bar{\zeta}_\theta)}{|\Psi_\zeta||\zeta_\theta|^3} + \frac{\Image(\Psi_{\zeta\zeta}\bar{\Psi}_\zeta \zeta_\theta)}{|\Psi_\zeta|^3|\zeta_\theta|}.
\]
\Add{Using the expression in \cref{eq:taylor-psi} and }taking $\zeta = O(\eps)$ into account, each component can be calculated as follows:
\begin{align*}
    |\Psi_\zeta|^2 &= |1 + a\zeta + b\zeta^2/2 + O(\zeta^3)|^2 \\
    &= 1 + 2\Real(a\zeta + b\zeta^2/2) + |a\zeta|^2 + O(\eps^3), \\
    |\zeta_\theta| &= \sqrt{1 + \eps^4 R_\theta^2}|\zeta|, \\
    \Image(\zeta_{\theta\theta}\bar{\zeta}_\theta) &= \Image\left((\eps^4 R_\theta^2 + \eps^2 R_{\theta\theta} - 1 + 2i\eps^2 R_\theta)(\eps^2 R_\theta - i)\right)|\zeta|^2 \\
    &= (1 - \eps^2 R_{\theta\theta} + O(\eps^4))|\zeta|^2, \\
    \Image(\Psi_{\zeta\zeta}\bar{\Psi}_\zeta \zeta_\theta) &= \Image\left((a + b\zeta + O(\zeta^2))\overline{(1 + a\zeta + b\zeta^2/2 + O(\zeta^3))}(\eps^2 R_\theta + i)\zeta\right) \\
    &= \Real(a\zeta) + |a\zeta|^2 + \Real(b\zeta^2) + O(\eps^3).
\end{align*}
We then obtain that
\begin{align*}
    \frac{\Image(\zeta_{\theta\theta}\bar{\zeta}_\theta)}{|\zeta_\theta|^3|\Psi_\zeta|} 
    &= \frac{1}{\eps}\left[1 \!-\! \Real(a\zeta) \!-\! \eps^2 R_{\theta\theta} \!-\! \frac{1}{2}(\Real(b\zeta^2) \!+\! |a\zeta|^2 \!-\! 3\Real(a\zeta)^2) \!-\! \eps^2 R \!+\! O(\eps^3)\right], \\
    \frac{\Image(\Psi_{\zeta\zeta}\bar{\Psi}_\zeta \zeta_\theta)}{|\zeta_\theta||\Psi_\zeta|^3} &= \frac{1}{\eps}\big[\Real(a\zeta) + |a\zeta|^2 + \Real(b\zeta^2) - 3\Real(a\zeta)^2 + O(\eps^3)\big].
\end{align*}
Therefore, we have
\begin{align*}
    \kappa &= \frac{1}{\eps}\bigg[1 - \eps^2(R_{\theta\theta} + R) + \frac{1}{2}\big(|a\zeta|^2 + \Real(b\zeta^2) - 3\Real(a\zeta)^2\big) + O(\eps^3)\bigg] \\
    &= \frac{1}{\eps}\bigg[1 - \eps^2(R_{\theta\theta} + R) + \frac{\eps^2}{2}\big(|a|^2 + \Real(be^{2i\theta}) - 3\Real(ae^{i\theta})^2\big) + O(\eps^3)\bigg].
\end{align*}
Now we use the facts
\begin{equation*}
    \sin(2\theta) = \sum_{k=0}^{\infty}\frac{8}{\pi(3 - 4k -4k^2)}\cos((2k + 1)\theta)\Erase{, \qquad b_2 - 3a_1 a_2 = K_s,}
\end{equation*}
and
\begin{equation*}
    b_2 - 3a_1 a_2 = K_s
\end{equation*}
\Add{from the definitions in \cref{eq:psi-ab}. }Then we obtain that
\begin{equation}
    \kappa = \frac{1}{\eps} - \eps\left(R_0 + \sum_{n \geq 2}(1 - n^2)R_n\cos(n\theta)\right) - \frac{\eps}{2}\sum_{n\geq 0}e_n(s)\cos(n\theta) + O(\eps^2), \label{kappa}
\end{equation}
where
\begin{equation*}
    e_n(s) = \begin{cases}
        \dfrac{1}{2}(a_1^2(s) + a_2^2(s)), & n = 0, \\[1em]
        -b_1(s) + \dfrac{3}{2}(a_1^2(s) - a_2^2(s)), & n = 2, \\[1em]
        \dfrac{8K_s(s)}{\pi(3 - 4k - 4k^2)}, & n = 2k + 1, k \geq 0, \\[1em]
        0 & n = 2k, k \geq 2.
    \end{cases}
\end{equation*}

This expression implies that
\begin{equation}
    \average{\kappa} = \frac{1}{\eps} - \eps R_0 - \frac{\eps}{4}(a_1^2 + a_2^2) + O(\eps^2), \label{average-kappa}
\end{equation}
that is, $\average{\kappa}$ extracts $0$-mode from $\kappa$.

\subsection{Motion along the boundary}
Substituting \cref{traveling-wave-speed,normal-velocity,kappa,average-kappa} into \cref{slow-interface-equation}, we obtain
\begin{equation*}
    \begin{aligned}
        &\dot{s}(\cos\theta + \eps B_1^\eps) + \eps^3 \dot{R}(1 + \eps B_2^\eps) \\
        =& \eps c'(v_c)w_\eps(R_0) + \Erase{\eps^2 C_1^\eps}\Add{\eps^{8}C_1^\eps} + \eps^2\sum_{n\geq 2}(1 - n^2)R_n\cos(n\theta) + \frac{\eps^2}{2}\sum_{n \geq 1}e_n\cos(n\theta) + \eps^3B_3^\eps,
    \end{aligned}
\end{equation*}
where $C_1^\eps$ represents the $0$-mode terms of order \Erase{$w_\eps^2$}\Add{\( \eps^0\)}, and $B_3^\eps$ collects the terms without the $0$-mode.
Dividing both sides by $1 + \eps B_2^\eps$ yields
\begin{equation}
    \begin{aligned}
        &\eps\left(\eps^2 \dot{R}_0 - c'(v_c)w_\eps(R_0)\right) + \left(\dot{s} - \frac{4\eps^2}{3\pi}K_s\right)\cos\theta + \eps^2\sum_{n \geq 2}E_n(s,t)\cos(n\theta) \\
        =& \Erase{\eps^2 C_2^\eps}\Add{\eps^{8}C_2^\eps} + \eps B_5^\eps \dot{s} + \eps^3 B_6^\eps,
    \end{aligned}
\end{equation}
where $C_2^\eps$ denotes the $0$-mode terms of order \Erase{$w_\eps^2$}\Add{\( \eps^0 \)}, $B_5^\eps$ and $B_6^\eps$ are composed of terms without a $0$-mode,
and $E_n(s,t)$ is defined by
\begin{gather*}
    E_n(s,t) = \eps \dot{R}_n(t) + (n^2 - 1)R_n(t) - \frac{e_n(s)}{2},
\end{gather*}
\Add{all of which arise due to the division by \( 1 + \eps B_2^\eps \).}
Let $B_{i,n}^\eps$ denote the $n$-th Fourier coefficient of $B_{i}^\eps$ for $i = 5,6$.
By projecting the equation onto each Fourier mode, we conclude that
\begin{equation}
    \begin{cases}
        \dot{s} = \dfrac{4\eps^2}{3\pi} K_s + \eps B_{5,1}^\eps \dot{s} + \eps^3 B_{6,1}^\eps, \\[0.85em]
        \eps^{2} \dot{R}_0 = c'(v_c)w_\eps(R_0) + \eps^{\Add{7}} C_2^\eps, \\[0.85em]
        \eps \dot{R}_n = -(n^2 - 1)R_n + \dfrac{e_n}{2} + \eps^{-1} B_{5,n}^\eps \dot{s} + \eps B_{6,n}^\eps, \quad n \geq 2.
    \end{cases} \label{eq-s,R}
\end{equation}
For the first equation, we solve for $\dot{s}$ by collecting the terms and dividing by $1 - \eps B_{5,1}^\eps$. This reveals that $\dot{s} = O(\eps^2)$. Furthermore, \Erase{assuming that $w_\eps$ is of order $\eps^l$ with $l \geq 0$, and} setting $\Erase{w}\Add{\tilde{w}_\eps}(R_0) = \eps^{-3}w_\eps(R_0) \Add{= O(1)}$, we arrive at the following system:
\begin{equation}
    \begin{cases}
        \dot{s} = \dfrac{4\eps^2}{3\pi} K_s + \eps^3 \tilde{B}_{6,1}^\eps, \\[0.85em]
        \eps^{\Add{-1}}\dot{R}_0 = c'(v_c)\Erase{w}\Add{\tilde{w}_\eps}(R_0) + \Erase{\eps^{1-l}}\Add{\eps^{4}} C_2^\eps, \\[0.85em]
        \eps\dot{R}_n = -(n^2 - 1)R_n + \dfrac{e_n}{2} + \eps \tilde{B}_{6,n}^\eps, \quad n \geq 2,
    \end{cases} \label{mab}
\end{equation}
where $\tilde{B}_{6,n}^\eps$ represents the remainder terms. Specifically, for $n=1$, $\tilde{B}_{6,1}^\eps$ accounts for the higher-order terms generated by the division by $1 - \eps B_{5,1}^\eps$. For $n \geq 2$, $\tilde{B}_{6,n}^\eps$ incorporates the corrections resulting from the substitution of $\dot{s}$, which is of order $\eps^2$. 
The first equation of the system \cref{mab} indicates that the leading-order dynamics of the interface position $s(t)$ are governed by
\[
\dot{s} = \frac{4\eps^2}{3\pi} K_s.
\]
This relationship provides a rigorous mathematical basis for the ``shape-sensing''. It implies that the interface is not passive but actively seeks out geometric features; specifically, the polarized domain is attracted to protrusions (local curvature maxima) and repelled from indentations. This mechanism allows the cell to orient itself purely based on the geometry of its environment, without external chemical gradients.

Finally, we investigate the stability of the equilibria for the problem \cref{mab}.
Letting $R = (R_0, R_2, R_3, \dots)$ by abuse of notation, we establish the following theorem concerning the existence and stability of the stationary solutions to \cref{mab}.
\begin{theorem*}
    Let $s^0$ be a non-degenerate critical point of the boundary curvature $K(s)$, i.e., $K_s(s^0) = 0$ and $K_{ss}(s^0) \neq 0$.
    Then, for sufficiently small $\eps > 0$, there exists a unique equilibrium solution $(s^\eps, R^\eps)$ to the system \cref{mab} in a neighborhood of the unperturbed equilibrium $(s^0, R^0)$.
    Moreover, the linear stability of $(s^\eps, R^\eps)$ is determined by the local geometry of the boundary at $s^0$:
    \begin{enumerate}
        \item If $s^0$ is a local maximum of $K(s)$ (i.e., $K_{ss}(s^0) < 0$), the equilibrium is asymptotically stable.
        \item If $s^0$ is a local minimum of $K(s)$ (i.e., $K_{ss}(s^0) > 0$), the equilibrium is unstable.
    \end{enumerate}
\end{theorem*}

\begin{proof}
Let $P(s, R, \eps)$ denote the vector field given by the right-hand side of the system \cref{mab}, rescaled such that the leading-order terms are of $O(1)$.
Namely, $P(s, R, \eps) = (P_0, P_1, P_2, \dots)$ where
\begin{align*}
    P_0(s, R, \eps) &= \frac{4}{3\pi}K_s(s) + \eps \tilde{B}_{6,1}^\eps, \\
    P_1(s, R, \eps) &= c'(v_c)\Erase{w}\Add{\tilde{w}_\eps}(R_0) + \Erase{\eps^{1-l}}\Add{\eps^{4}} C_2^\eps, \\
    P_2(s, R, \eps) &= -(2^2 - 1)R_2 + \frac{e_2}{2} + \eps \tilde{B}_{6,2}^\eps, \\
    &\vdots \\
    P_n(s, R, \eps) &= -(n^2 - 1)R_n + \frac{e_n}{2} + \eps \tilde{B}_{6,n}^\eps, \\
    &\vdots
\end{align*}
An equilibrium $(s^0, R^0)$ of the unperturbed system ($\eps = 0$) is a solution to $P(s, R, 0) = 0$,
which must satisfy
\begin{gather*}
    K_s(s^0) = 0, \\
    \Erase{w}\Add{\tilde{w}_\eps}(R_0^0) = 0, \text{ i.e., } v_0(R_0^0) = v_c, \\
    R_n^0 = \frac{e_n(s^0)}{2(n^2 - 1)} \Erase{= 0} \text{ for } n \geq 2.
\end{gather*}
We now examine the unique existence and stability of an equilibrium for the perturbed system ($\eps > 0$) near $(s^0, R^0)$.
The Jacobian matrix of $P$ with respect to $(s, R)$ evaluated at $(s^0, R^0, 0)$ forms a lower triangular matrix:
\begin{equation*}
    \begin{pmatrix}
        \dfrac{4}{3\pi}K_{ss}(s^0) & 0 & 0 & \cdots & 0 & \cdots  \\[1em]
        0 & c'(v_c)\Erase{w'}\Add{\tilde{w}_\eps'}(R_0^0) & 0 & \cdots & 0 & \cdots \\[1em]
        \dfrac{e_2'(s^0)}{2} & 0 & -3 & \cdots & 0 & \cdots \\[1em]
        \vdots & \vdots & \vdots & \ddots & \vdots & \ddots \\[1em]
        \dfrac{e_n'(s^0)}{2} & 0 & 0 & \cdots & -(n^2 - 1) & \cdots \\[1em]
        \vdots & \vdots & \vdots & \ddots & \vdots & \ddots \\[1em]
    \end{pmatrix}.
\end{equation*}
The eigenvalues of the above matrix are given by the diagonal entries: $\frac{4}{3\pi}K_{ss}(s^0)$, $c'(v_c)\Erase{w'}\Add{\tilde{w}_\eps'}(R_0^0)$, and $-(n^2 - 1)$ for $n\geq 2$.
Since we have $c'(v_c)\Erase{w'}\Add{\tilde{w}_\eps'}(R_0^0) < 0$ \Add{at \( \eps = 0 \)} (see the \Erase{defenition}\Add{definition} of $\Erase{w}\Add{\tilde{w}_\eps}$ and Subsection \ref{condition for w}), the Implicit Function Theorem guarantees
that if $K_{ss}(s^0) \neq 0$, a unique solution to $P(s, R, \eps) = 0$ exists in a neighborhood of $(s,R) = (s^0, \Erase{R_0^0}\Add{R^0})$ for sufficiently small $\eps > 0$.
Furthermore, the stability of this solution for $\eps > 0$ is the same as that of $(s^0, R^0)$ and is therefore determined by the
sign of $K_{ss}(s^0)$.
Namely, the equilibrium for $\eps > 0$ is stable if $K_{ss}(s^0) < 0$, and is unstable if $K_{ss}(s^0) > 0$.
\end{proof}

\section{Concluding Remarks}

In this paper, we have derived a free boundary problem approximating the front behavior of the wave-pinning model valid up to the $O(\varepsilon^{-1})$ timescale via matched asymptotics.
Our results analytically corroborate the emergence of area-preserving mean curvature flow on the $O(\varepsilon^{-1})$ timescale, a behavior that has been suggested in previous studies.
Furthermore, we demonstrated that the interface undergoes a slow drift along the domain boundary on the $O(\varepsilon^{-2})$ timescale, driven by the interaction between curvature flow and boundary geometry.

While these findings provide theoretical insights into the long-time behavior of the wave-pinning model, several mathematical challenges remain.
A rigorous convergence analysis of the original system to the derived interface equations as $\varepsilon \downarrow 0$ is required.
For such rigorous treatments of reaction-diffusion systems and their interfaces, see \cite{gomez2023front} and references therein.
\Erase{In Section 3, we utilized a scaling argument to justify the assumption that the first-order correction $v_1$ is of order $\varepsilon^2$; however, verifying this assumption requires examining the properties of the transmission problem for $v_1$, which necessitates a more delicate argument.
Verifying this assumption would likely require establishing $L^\infty$ estimates for the transmission problem on the perturbed domain, which involves subtle regularity issues at the boundary contact points.
While we leave this rigorous proof for future work, the consistency of our asymptotic results with numerical simulations strongly supports the validity of this scaling ansatz.}
\Add{In Section 3, the relevant quantity is not the absolute size of \( v_1 \), whose constant mode is fixed by the mass constraint \cref{fbp-mc1}, but rather the oscillatory part \( v_1 - \langle v_1 \rangle \) on \( \Gamma \).
We used a formal scaling argument to absorb this term into the \( O(\eps^3) \) remainder. A rigorous justification would require establishing gradient estimates for the associated elliptic transmission problem on the perturbed domain, including a careful treatment of the boundary contact points.
We leave this problem for future work.}
In addition to these analytical challenges, the numerical treatment of the $v_1$ term presents significant difficulties.

These complexities associated with $v_1$ lead us to consider a simplified model in which the $v_1$ term is neglected.
This simplified model retains a sufficient mechanism for wave-pinning.
Consequently, this simplified model serves not only as a valid approximation of the full wave-pinning model but also as a minimal framework that distills its essence, providing a tractable mathematical description suitable for investigating long-time behaviors.
Detailed analysis of this simplified model also remains an interesting subject for future research.

We discuss the implications of our results for understanding cell motility models, specifically the model proposed by Camley et al. \cite{camley2017crawling}.
Their model, which couples the wave-pinning reaction-diffusion system with phase-field mechanics, exhibits a transition between straight and rotating motion depending on parameter values.
Camley et al. hypothesized that ``shape sensing'' inherent to the wave-pinning kinetics plays a crucial role in the mechanism of rotation.
Our derivation of the drift effect toward high-curvature regions provides analytical validation for this hypothesis.
In the context of rotating motion, the mechanism can be interpreted as a mismatch between mechanical and chemical stable states.
While mechanical forces driven by protrusion tend to deform the cell into a bilaterally symmetric shape (typically symmetric about the short axis), the internal reaction-diffusion dynamics drives the Rho GTPase interface toward regions of higher curvature (i.e., the tips of the long axis), as shown in our $O(\varepsilon^{-2})$ analysis.
This continuous breaking of symmetry likely induces the rotational behavior observed in numerical simulations.
Conversely, the existence of straight motion presents a mathematically intriguing phenomenon.
In a migrating cell, the polarity is typically localized at the leading edge, which corresponds to a region of relatively low curvature compared to the sides.
According to our analysis, such a configuration should be unstable, as the interface naturally drifts toward higher curvature.
The persistence of straight motion, therefore, suggests that the kinetics of the moving boundary (protrusion forces) exert a stabilizing effect that counteracts the geometric drift.
In other words, the coupling with boundary mechanics stabilizes a configuration that is inherently unstable within the pure wave-pinning framework.
Unraveling the mathematical structure of this competition—specifically, proving the existence and bifurcation of straight versus rotating solutions in the coupled system—remains a challenging and significant subject for future research.

Ultimately, our analysis suggests that cell polarity is a robust phenomenon where geometry acts as a governing cue. By reducing the complex reaction-diffusion system to a set of geometric evolution laws, we provide a theoretical tool that can be used to predict cell behavior in engineered micro-environments, bridging the gap between mathematical theory and experimental cell biology. \color{black}

\section*{Acknowledgments}
\Add{The authors are grateful to the anonymous referees for their thorough reading of the manuscript and for their insightful comments, which substantially contributed to improving the quality of this work.}
This work was supported by JSPS KAKENHI Grant Numbers JP22K03425, JP22K18677, JP23H00086, JP24K16964, 25K00918.


\end{document}